\documentclass{amsart}
\usepackage{graphicx}
\usepackage{caption}
\usepackage{subcaption}
\usepackage{float}
\usepackage{hyperref}
\hypersetup{
colorlinks=true,
linkcolor=red,
filecolor=magenta, 
urlcolor=cyan,
} 
\usepackage{color}

\usepackage{tikz}
\usetikzlibrary{shapes.geometric, arrows}
\tikzstyle{startstop} = [rectangle, rounded corners, minimum width=3cm, minimum height=1cm,text centered, draw=black, fill=red!30]
\tikzstyle{io} = [trapezium, trapezium left angle=70, trapezium right angle=110, minimum width=3cm, minimum height=1cm, text centered, draw=black, fill=blue!30]
\tikzstyle{decision} = [diamond, minimum width=3cm, minimum height=1cm, text centered, draw=black, fill=green!30]
\tikzstyle{process} = [rectangle, minimum width=3cm, minimum height=1cm, text centered, draw=black, fill=orange!30]
\tikzstyle{decision} = [diamond, minimum width=3cm, minimum height=1cm, text centered, draw=black, text width = 3cm, fill=green!30]
\tikzstyle{arrow} = [thick,->,>=stealth]

\usepackage {verbatim}
\usepackage{hyperref}
\usepackage[]{algorithm2e}
\usepackage{amsmath,amssymb,eucal,amsfonts,amsthm,enumerate}

\newcommand{\spacing}[1]{\renewcommand{\baselinestretch}{#1}\large\normalsize}
 \spacing{1.66} 
\newtheorem{thm}{Theorem}

\newtheorem{rem}{Remark}[section]

\newtheorem{Exp}{Example}

\newcommand{\R}{\mathbb{R}}

\newcommand{\beq}{\begin{equation}}
\newcommand{\eeq}{\end{equation}}
\newcommand{\grad}{\nabla}
\newcommand{\abs}[1]{\lvert#1\rvert}

\renewcommand{\O}{\Omega}

\newcommand{\pd}[2]{\frac{\partial #1}{\partial #2}}

\renewcommand{\L}{\mathcal{L}^2(\Omega)}
\newcommand{\Linfty}{\mathcal{L}^\infty(\Omega)}

\newcommand{\SWO}{\mathcal{W}^{1,2}_0(\Omega)}
\newcommand{\SW}{\mathcal{W}^{1,2}(\Omega)}

\newcommand{\Lnorm}[1]{\abs{\abs{#1}}_{\mathcal{L}^2}}
\newcommand{\Lonenorm}[1]{\abs{\abs{#1}}_{\mathcal{L}^1}}

\newcommand{\Lprod}[2]{(#1,#2)_{\mathcal{L}^2}}

\newcommand{\lgrad}[1]{\nabla_{\mathcal{L}^2}^{#1}}
\newcommand{\hgrad}[1]{\nabla_{\mathcal{H}^1}^{#1}}

\newcommand{\D}{\mathcal{D}}

    \usepackage{newunicodechar}
    \newunicodechar{ﬁ}{fi}
    \newunicodechar{ﬀ}{ff}

\begin{document}

\title{Parameter estimation in an elliptic problem.}

\author{Abinash Nayak}
\address{ Department of Mathematics, University of Alabama at Birmingham, Campbell Hall, Rm. 452, 1300 University Blvd., Birmingham, AL 35233}
%    Current address Department of Mathematics, University of Alabama at Birmingham, Birmingham, Alabama
%\curraddr{Department of Mathematics, UAB, Campbell Hall, Rm. 452, 1300 University Blvd., Birmingham, AL 35233}
\email{nash101@uab.edu; avinashnike01@gmail.com}
% \thanks{*I am very thankful to my advisor for his support, encouragement and expert guidance throughout my research. }

    % Information for second author
% \author{Author Two}
% \address{Mathematical Research Section, School of Mathematical Sciences,
% Australian National University, Canberra ACT 2601, Australia}
% \email{two@maths.univ.edu.au}
% \thanks{Support information for the second author.}

%    General info
% \subjclass[2000]{Primary 54C40, 14E20; Secondary 46E25, 20C20}
%\thanks{This work was supported by the National Science Foundattion.}
% \date{January 1, 2001 and, in revised form, June 22, 2001.}
% \today
\date{\today}
%\dedicatory{This paper is dedicated to my advisor, Dr. Ian 
%Knowles.}


\begin{abstract}
A new variational approach to solve the problem of estimating the (possibly discontinuous) coefficient functions $p$, $q$ and $f$ in elliptic equations of the form $-\nabla \cdot (p(x)\nabla u) + \lambda q(x) u = f$, $x \in \O \subset \R^n$, from a knowledge  of the solutions $u_\lambda$.
\end{abstract}

\maketitle
\begin{Exp}\label{Example Parameter identification}\textbf{Parameter identification.}\\
By parameter identiﬁcation one usually denotes the problem of reconstructing the unknown coefficients in a partial diﬀerential equation from (indirect) measurements of the solution or a noisy solution. A simple example is the following model from groundwater filtration, which is modeled through the following elliptic equation 
\begin{equation}\label{parameter identification eq.}
    -\grad.(p\grad u) = f, 
\end{equation}
in $\Omega \subset \R^d$, where $u$ is the unknown, $f$ a given source, and $p$ the hydraulic permittivity. The direct problem consists in solving the partial differential equation for $u$, given a $p$ and suitable boundary conditions on $\partial \O$. The inverse problem consists in reconstructing the unknown parameter $p$ on $\O$ given a noisy measurement of the solution,
\beq
    u_\delta(x) = u(x) + \epsilon_\delta(x), \hspace{1cm} x \in \O.
\eeq
If the solution of the direct problem is unique for each parameter $p$, which is the case for the groundwater filtration problem with appropriate boundary conditions, then one can introduce the \textit{parameter-to-solution map}, $p \rightarrow u_p$, where $u_p$ is the solution to the direct problem given a specific $p$. Note that even if the direct problem is linear (for $u$), the inverse problem and the parameter-to-output map are usually non-linear. For example, in the ground water filtration problem we have $u_{2p} = \frac{1}{2}u_p$, and not $u_{2p} = 2u_p$ and hence, the problem is not linear.

The uniqueness question for parameter identification problems is usually denoted as \textit{identifiability}. For example, if $\O = [0, 1]$ then integrating equation \eqref{parameter identification eq.} yields
\beq \label{par. idn. [0,1] eq.}
    -p(x)u'(x) + p(0)u'(0) = \int_0^x f(\xi) d\xi.
\eeq
Hence, from \eqref{par. idn. [0,1] eq.}, the parameter $p$ can be uniquely determined (a.e.) for a given $u$ and $f$ provided $u' \neq 0$ (a.e.) and knowing $p(0)$, see \cite{Knowles+Yan,Knowles(PIFEP),Knowles+Yan_1,Knowles+Yan+Le,Knowles+Yan_2} for inverse problems related to ground water modelling.

The naive approach to retrieve the parameter 
\beq
    p(x) = \frac{p(0)u'(0) - \int_0^xf(\xi)d\xi}{u'(x)}
\eeq
shows that besides the usual \textit{linear ill-posedness} arising from the fact that the data (usually noisy, $u_\delta$) have to differentiated, there is also a \textit{nonlinear ill-posedness} arising from the quotient, whose consequence is that errors at the small values of $u'$ are amplified much stronger than errors at large values of $u'$. That is, if $u'(x)$ is very small in an interval $I$ then, though we have identifiability, in practice we must expect very high error due to the noise amplification. 
\end{Exp}

\section{Parameters estimation for elliptic partial differential equation}\label{Parameter estimation for elliptic partial differential equation}

The problem of estimating (or identifying) the coefficients (or parameters) involved in a elliptic differential equation has a paramount practical significance. Mathematically, it can be considered as finding the parameters ($p$, $q$, $f$) of the following differential equation, for the known solution $u_\lambda$ (depending on a parameter $\lambda$), 
\beq\label{elliptic eq.}
    L_\lambda u \equiv - \nabla \cdot (p(x) \nabla u) + \lambda q(x) u = f(x), \;\; x \in \O
\eeq
where $\O$ is an open simply connected bounded set with a $C^1$-boundary in $\R^n$. Here we assume the parameters
\beq\label{p,q,f domains}
    f \in \L, \;\; q \in \Linfty, \;\; \mbox{and } p \in \Linfty
\eeq
with $p$ satisfying 
\beq\label{ellipticity cond.}
    p(x) \geq \nu > 0, \;\; x \in \O,
\eeq
and $q$, with the real parameter $\lambda$, is such that the homogeneous ($f \equiv 0$) Dirichlet operator $L_\lambda = L_{\lambda, p,q}$ (i.e., $L_\lambda$ acting on $\SWO$) satisfies\footnote{note that as $q \in \Linfty$, for $|\lambda|$ small enough, the condition \eqref{Llambda positive} is true.}
\beq\label{Llambda positive}
    L_\lambda \text{ is a positive operator } \L.
\eeq
It is known, in \cite{Gilbarg_Trudinger}, that for a $\phi \in \SW$ the generalized Dirichlet problem associated with \eqref{elliptic eq.}, with boundary condition 
\beq
    u|_{\partial \O} = \phi|_{\partial \O},  
\eeq
is uniquely solvable, and that the solution lie in the Sobolev space $\SW$. We are interested here in the corresponding inverse problem: given the solutions $u_\lambda$ (for one or more values of $\lambda$), find one or more of the coefficient functions $p$, $q$ and $f$. This inverse problem of identifying the parameters is an (non-linear) ill-posed problem, arising from the fact that it involves differentiation of (noisy) data.

Such inverse problems are of interest in connection with groundwater flow (and also oil reservoir simulation); see \cite{Knowles(PIFEP), Knowles+Yan, Knowles+Yan_2, Knowles+Yan_1, Knowles+Yan+Le} and the references therein. In such cases the flow in the porous medium is governed by the following diffusion equation
\beq\label{groundwater model}
 \nabla \cdot (P(x) \nabla w(x,t)) = S(x) \pd{w}{t} - R(x,t),
\eeq
in which $w$ represents the piezometric head, $P$ the hydraulic conductivity (or sometimes, for a two-dimensional aquifer, the transmissivity), $R$ the recharge, and $S$ the storativity of the aquifer. In the case that the aquifer reaches a steady-state condition, we have that $\pd{w}{t} = 0$ and $R=R(x)$, which is essentially the equation in \eqref{elliptic eq.}.

In \cite{Knowles1999} the theoretical framework was given for a general approach to the problem of computing, from a knowledge of the piezometric head values $w(x,t)$ of the aquifer over space and time, reliable values for the aquifer parameters. The basic idea in \cite{Knowles1999} is to transform, by appropriate means\footnote{for example, Finite Laplace transform on the time $t$ variable}, data from solutions of \eqref{groundwater model} to solution values $u_\lambda(x)$ of the following elliptic equation
\beq
    - \nabla \cdot (P(x)\nabla u) + \lambda S(x) u = F(x,\lambda), \;\; x \in \O,
\eeq
where $\lambda$ is a transform parameter and $F$ depends on $R$, $S$ and $\lambda$ in a known way. The triplet $(P,S,F)$ is then found (under suitable conditions on the solutions $u_\lambda$ and the form of  $R$) as the unique global minimum of a certain convex functional, which is discussed below.

These parameters are estimated by minimizing a (strictly) convex functional, whose minimizers corresponds to the original parameter triplet $(P,Q,F)$. The functional used in \cite{Knowles(PIFEP), Knowles+Yan, Knowles+Yan_2, Knowles+Yan_1, Knowles+Yan+Le, Knowles1999} can be generalized as follows: let a solution(s) $u_\lambda$ (depending on $\lambda$) of \eqref{elliptic eq.} be known (given) for which $(P,Q,F)$ are the coefficients corresponding to $p$, $q$ and $f$, respectively, that we seek to recover. For any $c = (p,q,f)$, where $p$, $q$ and $f$ satisfying \eqref{p,q,f domains}, \eqref{ellipticity cond.} and \eqref{Llambda positive}, let $u_{\lambda,c}$ (depending on $\lambda$) denote the solution of \eqref{elliptic eq.} corresponding to the choice of $c = (p,q,f)$ with the boundary condition 
\beq
    u_{\lambda,c}|_{\partial \O} = u_\lambda|_{\partial \O}.
\eeq
Thus, we have $u_\lambda = u_{\lambda,\hat{c}}$, for $\hat{c} = (P,Q,F)$. It's proved in \cite{Knowles2000} that $\hat{c} = (P,Q,F)$ is a minimizer of the following convex functional 
\begin{align}
    G_\lambda(c) &= \Lprod{L_{\lambda,p,q(u_\lambda - u_{\lambda,c})}}{u_\lambda - u_{\lambda,c}} \\
    &= \int_\O p(x)|\nabla(u_\lambda - u_{\lambda,c})|^2 + \lambda q(x)(u_\lambda- u_{\lambda,c})^2 \; dx,
\end{align}
where $c \in \D(G_\lambda) := \{ (p,q,f)|\; p,q,f \text{ satisfy } \eqref{p,q,f domains}, \eqref{ellipticity cond.}, \eqref{Llambda positive} \text{ and } p|\Gamma = P|_\Gamma \}$, where $\Gamma$ is a hypersurface in $\O$ transversal to $\nabla u_\lambda$. It is convenient to take $\Gamma$ to be the boundary of the bounded region $\O$, and henceforth we assume this to be so. We state some of the properties of the functional $G$, from \cite{Knowles2000},

\begin{thm}\mbox{ }\\
\begin{enumerate}
    \item For any $c = (p,q,f) \in \D(G_\lambda)$,
    \begin{gather}
        G_\lambda(x) = \int_\O p(x)(|\nabla u_\lambda|^2 - |\nabla u_{\lambda,c}|^2) + \lambda q(x)(u_\lambda^2 - u_{\lambda,c}^2) \notag \\
        - 2f(x)(u_\lambda - u_{\lambda,c}) \;dx.
    \end{gather}
    
    \item $G_\lambda(c) \geq 0$ for all $c \in \D(G_\lambda)$, and $G(c) = 0$ if and only if $u_\lambda = u_{\lambda,c}$.
    
    \item For $c_1 = (p_1,q_1,f_1)$ and $c_2 = (p_2,q_2,f_2)$ in $\D(G_\lambda)$, we have 
    \begin{align}
        G_\lambda(c_1) - G_\lambda(c_2) &= \int_\O (p_1 - p_2)(|\nabla u_\lambda|^2 - \nabla u_{\lambda,c_1} \cdot \nabla u_{\lambda, c_2}) + \lambda(q_1 - q_2)(u_\lambda^2 - \nonumber \\
        & \hspace{1cm} u_{\lambda,c_1}u_{\lambda,c_2}) - 2(f_1 - f_2)(u - \frac{u_{c_1} + u_{c_2}}{2}) \; dx
    \end{align}{}
    
    \item The first G$\hat{a}$teaux differential\footnote{can be proved that it is also the first Fr$\acute{e}$chet derivative of $G_\lambda$ at $c$} for $G_\lambda$ at any $c \in \D(G_\lambda)$ is given by 
    \begin{align}\label{Gp elliptic}
        G_\lambda'(c) [h_1, h_2, h_3] =& \int_\O (|\nabla u_\lambda|^2 - |\nabla u_{\lambda, c}|^2)h_1 + \lambda (u_\lambda^2 - u_{\lambda,c}^2)h_2 \nonumber \\
        & \hspace{2cm} -2(u_\lambda - u_{\lambda,c})h_3 \; dx,
    \end{align}{}
    for $h_1, h_2 \in \Linfty$ with $h_1|_{\partial \O} = 0$, and $h_3 \in \L$, and $G_\lambda'(c) = 0$ if and only if $u_\lambda = u_{\lambda,c}$.
    
    \item The second G$\hat{a}$teaux differential\footnote{can be proved that it is also the second Fr$\acute{e}$chet derivative of $G_\lambda$ at $c$} of $G_\lambda$ at any $c \in \D(G_\lambda)$  is given by
    \beq \label{Gpp elliptic}
        G''(c)[h,k] = 2\Lprod{L_{\lambda,p,q}^{-1}(e(h))}{e(k)},
    \eeq
    where $h = (h_1,h_2,h_3)$, $k = (k_1,k_2,k_3)$, and the functions $h_1,h_2, k_1, k_2 \in \Linfty$, with $h_1|{\partial \O} = k_1|_{\partial \O} = 0$, $h_3, k_3 \in \L$, and 
    \beq
        e(h) = - \nabla \cdot (h_1 \nabla u_{\lambda,c}) + \lambda h_2 u_{\lambda,c} - h_3.
    \eeq
\end{enumerate}{}
\end{thm}{}

\begin{rem}
The convexity of the functional $G_\lambda$ can be seen from \eqref{Gpp elliptic}, as for $h = k$ we have 
\beq \label{Gpp elliptic 1}
    G''(c)[h,h] = 2\Lprod{L_{\lambda,p,q}^{-1}(e(h))}{e(h)},
\eeq
and by the positivity of $L_{\lambda,p,q}$ for any $c = (p,q,f) \in \D(G_\lambda
)$ we get $G_\lambda''(c)[h,h] \geq 0$, but this doesn't imply the strict convexity as we can have $e(h) = 0$ for $h \not\equiv 0$, i.e., not all $h_1$, $h_2$ and $h_3$ are zeros simultaneously. It does make sense, as one can not expect to inversely recover three (unknown) parameters $(P,Q,F)$ through solving only one equation \eqref{elliptic eq.}, for a particular $u_\lambda$. However, as proved in \cite{Knowles_LaRussa}, if one has solutions $u_\lambda$'s corresponding to certain $\lambda$'s $\in I$ (an index set), then one can have a combination of the convex $G_\lambda$'s to obtain a strictly convex functional $G$, i.e., 
\beq
    G = \sum_{\lambda \in I} G_\lambda.
\eeq
Intuitively (clearly for $\O \subset \R$), one can see that there need to be at least three $\lambda$'s in $I$ such that the following system of equations has a unique solution $(P,Q,F)$, for known $u_\lambda$'s,
\begin{align*}
    -\nabla \cdot (p(x) \nabla u_{\lambda_1}) + \lambda_1 q(x) u_{\lambda_1} &= f(x)\\
    -\nabla \cdot (p(x) \nabla u_{\lambda_2}) + \lambda_2 q(x) u_{\lambda_2} &= f(x)\\
    -\nabla \cdot (p(x) \nabla u_{\lambda_3}) + \lambda_3 q(x) u_{\lambda_3} &= f(x),
\end{align*}{}
for $x \in \O$. If the above holds, then we have (by linearity) $G''(c)[h,h] = \sum_{\lambda \in I} G_\lambda''(c)[h,h] \geq 0$, however, here $G''(c)[h,h] = 0$ if and only if $G_{\lambda}''(c)[h,h] = 0$ for all $\lambda$'s in $I$ and hence, $h_1 = h_2 = h_3 \equiv 0$, i.e., $G$ is strictly convex, for details (especially, when $\O \subset \R^n$, $n > 1$) see \cite{Knowles_LaRussa}.  
\end{rem}{}

\begin{rem}
Similar to the analysis performed for the previous regularization methods, one can observe that, through a descent algorithm, there exists a sequence of functions $c_m = (p_m, q_m, f_m) \in \D(G)$ such that $c_m$ converges to the unique (global) minimum $(P,Q,F)$ of the functional $G$, i.e., $p_m$, $q_m$ and $f_m$ converges weakly to $P$, $Q$ and $F$ in $\Linfty$, $\Linfty$ and $\L$, respectively.
\end{rem}{}

\begin{rem}
Note that during the descent method one needs to preserve the (given) boundary information of the parameter $P$, i.e., for the recovery of parameter $P$ the boundary data $p_m|_{\Gamma} = P|_{\Gamma}$, for all $m$, where $\Gamma = \partial \O$ should be invariant during the descent process. Hence, it leads to a constraint minimization, the constraint being to preserve $p_m|_{\partial \O} = P|_{\partial \O}$, for all $m$. This is achieved by having a gradient (descent) direction which vanishes at the boundary, which is given by the Neubereger gradient (or Sobolev gradient), see \cite{Neuberger1997}, chosen so that, 
\beq
    G_\lambda'(c_m)[h] = \Lprod{\hgrad{c_m}G_\lambda}{h}
\eeq
for all $h = [h_1,0,0]$, where $h_1 \in \SWO \cap \Linfty$. It can be computed by solving the following Dirichlet differential equation 
\begin{gather}\label{elliptic Dir. Neuberger}
    -\Delta g + g = \lgrad{c_m}G_\lambda \notag\\
    g|_{\partial \O} = 0,
\end{gather}
where $\lgrad{c_m}G_\lambda = |\nabla u_\lambda|^2 - |\nabla u_{\lambda,c_m}|^2$, from \eqref{Gp elliptic} and hence, we have the Neuberger or Sobolev gradient defined as
\beq
    \hgrad{}G_\lambda := g = (I - \Delta)^{-1}(\lgrad{}G_\lambda),
\eeq
which is zero on the boundary $\partial \O$. Thus, during the descent process the sequence $p_{m+1} = p_m - \alpha \hgrad{c_m}G_\lambda$, for an appropriate $\alpha$, converges weakly to $P$ in $\L$, with $p_m|_{\partial \O} = p_0|_{\partial \O} = P|_{\partial \O}$ (invariant).
\end{rem}{}

\begin{rem}
During the descent algorithm, in this scenario, one follows the directional descends at a particular stage, i.e, at $c_m = (p_m,q_m,f_m)$ the functional $G$ is first minimized in any particular direction (say $h = (h_1,0,0)$) to get $p_{m+1} = p_m - \alpha \hgrad{c_m}G$, for an appropriate $\alpha$, where the directional derivative $\hgrad{c_m}G = (I - \Delta)^{-1}(\lgrad{c_m}G)$ and $\lgrad{c_m}G = \sum_{\lambda \in I} (|\nabla u_\lambda|^2 - |\nabla u_{\lambda,c_m}|^2)$, and $c_m$ is partially developed to $c_m^{(1)} = (p_{m+1}, q_m, f_m)$; then $G$ is minimized in another direction (say $h = (0,h_2,0)$) to get $q_{m+1} = q_m - \alpha \lgrad{c_m^{(1)}}G$, for an appropriate $\alpha$ and $\lgrad{c_m^{(1)}}G = \sum_{\lambda \in I} \lambda (u_{\lambda}^2 - u_{\lambda, c_m^{(1)}}^2)$ and again, $c_m^{(1)}$ is further improved to $c_m^{(2)} = (p_{m+1}, q_{m+1}, f_m)$; and then $G$ is minimzed in the last direction $h = (0,0,h_3)$ to get $f_{m+1} = f_m - \alpha \lgrad{c_m^{(2)}}G$, for an appropriate $\alpha$ and $\lgrad{c_m^{(2)}}G = \sum_{\lambda \in I} -2(u_\lambda - u_{\lambda, c_m^{(2)}})$, and finally, $c_m^{(2)}$ is updated for the next iteration to $c_{m+1} = (p_{m+1}, q_{m+1}, f_{m+1})$.
\end{rem}{}

\begin{rem}
This regularization method requires not only the strict positivity of the parameter $p$, i.e., $0 < \nu \leq p(x)$ for all $x \in \O$, but also an a-priori knowledge on the lower bound. This is very essential when applying this regularization method numerically, as the during the descent algorithm the sequence $p_m$ (usually) tends towards zero and if not bounded away from zero, by a positive constant, it leads to instability and as a result blows up the numerical solver, see \cite{Knowles(PIFEP)}.  
\end{rem}{}

\subsection{Using the new regularization method}
One can implement the regularization method developed above in Chapter $??$ instead, with the linear operator being, for any $c = (p,q,f)$ satisfying \eqref{p,q,f domains} and $p|_{\partial \O} = P|_{\partial \O}$ (given),
\beq
    T_{\lambda,u_\lambda}(c) := L_{\lambda,p,q}(u_\lambda) - f,
\eeq
where $\lambda$ and $u_\lambda$ is given, and $L_{\lambda,p,q}(u_\lambda) = -\nabla \cdot (p\nabla u_\lambda) + \lambda q u_\lambda$. Hence the operator equation $??$ can be formulated as
\beq\label{elliptic op. eq.}
    T_{\lambda, u_\lambda} (c) = 0,
\eeq
for $c = (p,q,f)$ satisfying \eqref{p,q,f domains} and $p|_{\partial \O} = P|_{\partial \O}$, and the respective inverse problem as: find the $c$ satisfying \eqref{elliptic op. eq.}, given $\lambda$ and $u_\lambda$. The corresponding minimizing functional here is, for any $c = (p,q,f)$ satisfying \eqref{p,q,f domains} and $p|_{\partial \O} = P|_{\partial \O}$,
\beq\label{GT}
    G_T(c) = \Lnorm{T_{\lambda,u_\lambda}(c)}^2 + \Lnorm{(|\nabla v_{\lambda,c}|)}^2, 
\eeq
where $v_\lambda$ is the solution of the following Dirichlet problem
\begin{gather}
     -\Delta v_{\lambda, c} = T_{\lambda, u_\lambda}(c)\\
     v_{\lambda, c}|_{\partial \O} = 0. \notag
\end{gather}{}

\begin{rem}
Note that this scenario is different from the previous examples in the sense that, here for a noisy data $u_{\lambda,\delta}$ we do not have a noisy right hand side in the operator equation $??$, rather we have a perturbed operator $T_{\lambda,u_{\lambda,\delta}}$. Also note that the stability of the recovery, in this regularization method, depends on the $\Lnorm{T_{\lambda,u_{\lambda,\delta}} - T_{\lambda, u_\lambda}}$, which in return depends on the stable differention of the noisy data $u_{\lambda,\delta}$.
\end{rem}{}

\begin{rem}
Notice that for a noisy $u_{\lambda, \delta}$ one has to be very careful when computing $T_{\lambda, u_{\lambda, \delta}}(c)$, for any $c$, directly as it inherently contains the second derivative of the noisy data $u_{\lambda, \delta}$ and hence, would lead to serious noise amplifications. One can track the descend indirectly via $??$, which involves an integration of the operator $T$. The noisy effect can be further mitigated by having the Nuberger gradient instead of the $\L$-gradient during the descend process, as $\hgrad{}G = (1 - \Delta)^{-1}\lgrad{}G$, i.e., further smoothing. On the other hand, the function $v_{\lambda,c}$ is a much smoother function as $v_{\lambda,c} = -\Delta^{-1}(T_{\lambda,u_\lambda}(c)) = -\Delta^{-1}(-\nabla \cdot (p \nabla u_\lambda) + \lambda q u_\lambda - f)$ and, helps significantly in providing regularity to the inverse recovery.
\end{rem}{}

\begin{rem}
The greatest advantage of this regularization method over the earlier one is its independence on the knowledge of the lower bound for the parameter $P$. Unlike the previous method, where one needs to provide a lower (positive) bound for the parameter $P$ recovery, otherwise the solver crashes (see \cite{Knowles(PIFEP)}), here one does not need to specify any such bounds for any parameters recovery. In fact, as we see in Examples \ref{Example elliptic 2} and \ref{Example elliptic 3}, one can even recover the parameter $P$ having both positive and negative values, under certain constraints. However, if the parameter $P$ is zero over certain sub-domain $\O' \subset \O$, then it is not possible to recover $P$ uniquely.
\end{rem}{}

\section{Numerical Results}
For simplicity, we consider the inverse problem of recovering only one parameter in \eqref{elliptic eq.}. We focus on the inverse recovery of the parameter $P$, since it is the most difficult parameter to recover, as explained in Example \ref{Example Parameter identification}, and compare our results with the results obtained in \cite{Knowles(PIFEP)}. Since we are recovering only a single parameter, the inverse recovery is unique for a single solution $u = u_\lambda$ for any particular $\lambda$. As mentioned above, the ill-posedness in the problem is concentrated in the computation of $\nabla u_\delta$ from $u_\delta$. In consequence, the reliability and effectiveness of any proposed computational algorithm for this problem is directly dependent on how well the numerical differentiation is computed. Though in chapter $??$ we provide a very efficient method for computing numerical differentiation in one dimension, we did not get the time to extend it to its multi-dimension version. Hence, in the following examples, unless otherwise stated, either we assumed no error in the measured data, i.e., $\delta = 0$ (especially, when we are considering $P$ to have both negative and positive values, since, then the solution $u$ has singular values, see Example \ref{Example elliptic 4}), or, when $\delta \neq 0$, we fit a smooth surface (usually a degree five two-dimensional surface, using the curve fitting toolbox in MATLAB) through the noisy data $u_\delta$ and then differentiate the smooth surface $\tilde{u}$ to approximate $\nabla u$. In all of the examples we discretized the domain $\O = [-1,1]\times [-1,1]$ into evenly spread $49\times 49$ grid and used the MATLAB inbuilt PDE solvers for solving the PDEs. 

\begin{Exp}\label{Example elliptic 1}
In the first example we assumed the parameter $P$ defined as, in \cite{Knowles(PIFEP)},
\beq
    P_1(x,y) = 
    \begin{cases}
    2, & \mbox{ if } |x| < 0.5 \mbox{ and } |y| < 0.5 \\
    0, & \mbox{ otherwise},
    \end{cases}
\eeq
and the noisy free test data $u$ is constructed by solving \eqref{elliptic eq.} with the boundary function $\phi(x,y) = x + y + 4$ on $\partial \O$. We then contaminate the data $u$ with uniform error to get $u_\delta$ such that the relative error $\frac{\Lonenorm{u - u_\delta}}{\Lonenorm{u}} \approx 7\%$. Figure \ref{fig elliptic 1} shows the true $u$, the noisy $u_\delta$ and the smoothed $\tilde{u}$, which is obtained from fitting a degree five two-dimensional polynomial through $u_\delta$, respectively. Notice that the information present in true $u$ is completely lost in the presence of noise (which is kind of extreme in this case) and, though smoothing (with a degree five polynomial) lead to some resemblance with the true $u$, it misses the key features. Nevertheless, the recovery, as seen in Figure \ref{fig elliptic 1}, is still quite impressive.  The parameter $P$ is recovered through minimizing the functional $G_T$, as defined in \eqref{GT}, and using the Dirichlet Neuberger gradient $\hgrad{}G$, as defined in \eqref{elliptic Dir. Neuberger}. The relative error in the recovered $\tilde{P}$ is $\frac{\Lonenorm{P - \tilde{P}}}{\Lonenorm{P}} \approx 13.42 \%$. We compare our results with the results obtained in \cite{Knowles(PIFEP)}, which is shown in Figure \ref{knowles P recoveries}. 

\begin{figure}[ht]
    \centering
    \begin{subfigure}{0.4\textwidth}
        \includegraphics[width=\textwidth]{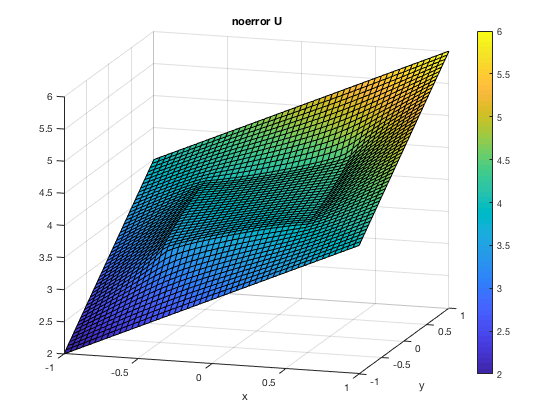}
        \caption{True $u$}
    \end{subfigure}
    \begin{subfigure}{0.4\textwidth}
        \includegraphics[width=\textwidth]{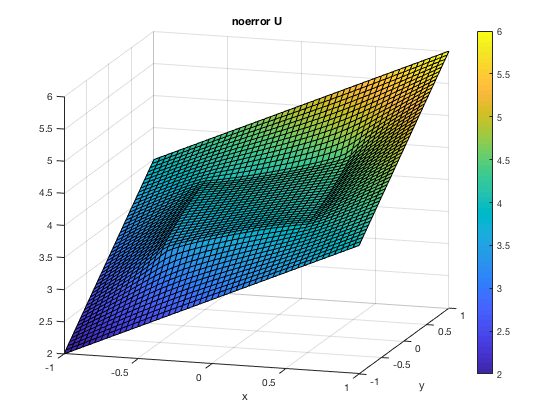}
        \caption{True $u$}
    \end{subfigure}

    \begin{subfigure}{0.4\textwidth}
        \includegraphics[width=\textwidth]{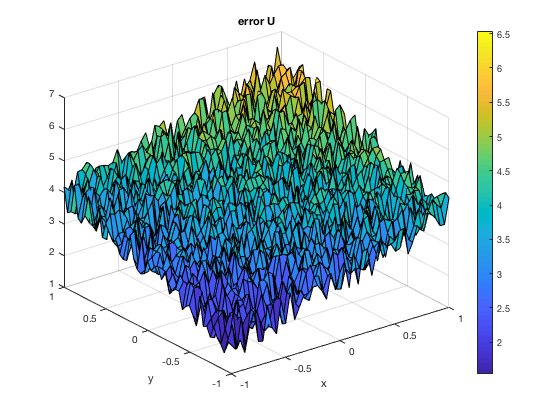}
        \caption{noisy $u_\delta$}
    \end{subfigure}
    \begin{subfigure}{0.4\textwidth}
        \includegraphics[width=\textwidth]{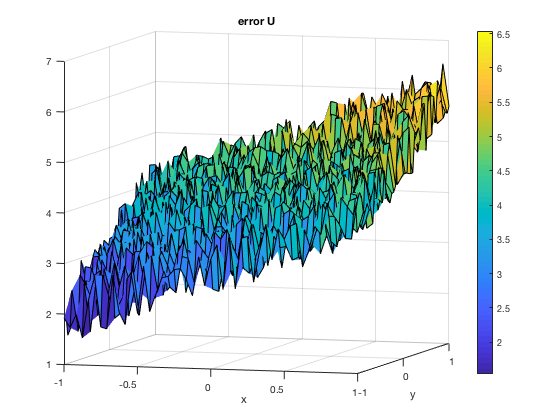}
        \caption{noisy $u_\delta$}
    \end{subfigure}

    \begin{subfigure}{0.4\textwidth}
        \includegraphics[width=\textwidth]{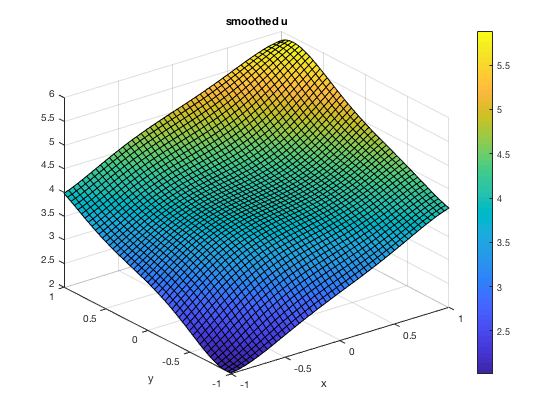}
        \caption{smoothed $u$}
    \end{subfigure}
    \begin{subfigure}{0.4\textwidth}
        \includegraphics[width=\textwidth]{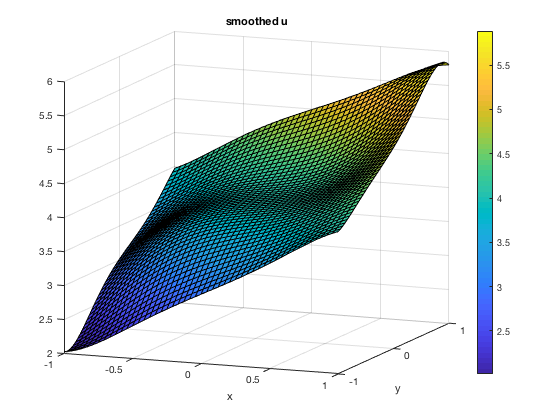}
        \caption{smoothed $u$}
    \end{subfigure}
    \caption{True $u$, noisy $u_\delta$ and smoothed $u$, for Example \ref{Example elliptic 1}.}  
    \label{fig elliptic 1}
\end{figure}

\begin{figure}[ht]
    \centering
    \begin{subfigure}{0.4\textwidth}
        \includegraphics[width=\textwidth]{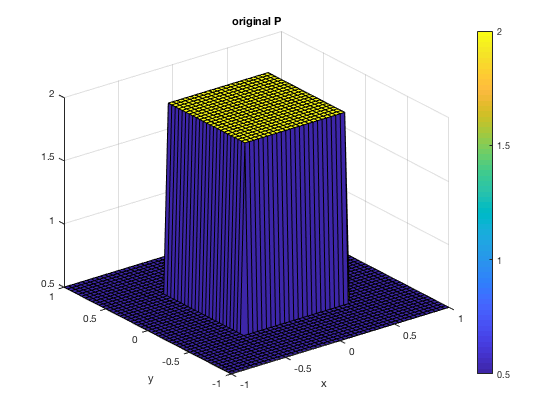}
        \caption{True $P$}
    \end{subfigure}
    \begin{subfigure}{0.4\textwidth}
        \includegraphics[width=\textwidth]{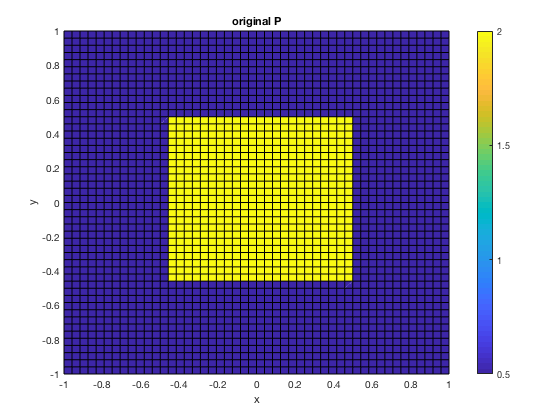}
        \caption{True $P$}
    \end{subfigure}

    \begin{subfigure}{0.4\textwidth}
        \includegraphics[width=\textwidth]{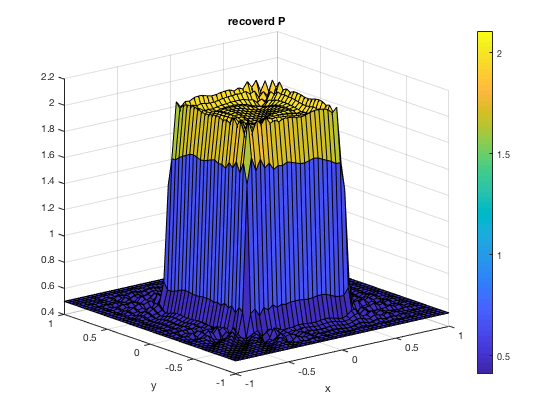}
        \caption{Recovered $\tilde{P}$}
    \end{subfigure}
    \begin{subfigure}{0.4\textwidth}
        \includegraphics[width=\textwidth]{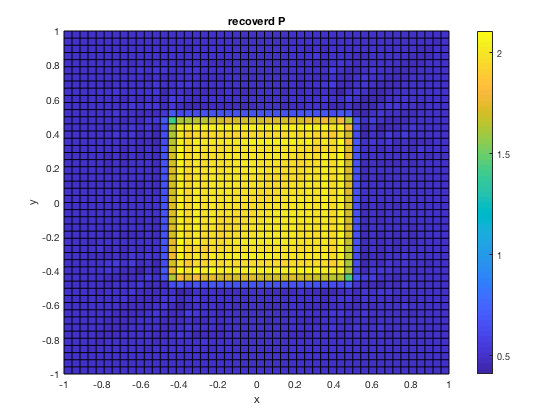}
        \caption{Recovered $\tilde{P}$}
    \end{subfigure}
    \caption{True $P$ and recovered $\tilde{P}$, for Example \ref{Example elliptic 1}.}  
    \label{true and recovered P}
\end{figure}

\begin{figure}[ht]
    \centering
    \begin{subfigure}{0.4\textwidth}
        \includegraphics[width=\textwidth]{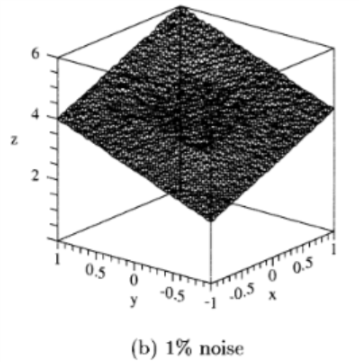}
        \caption{$1\%$ error in the solution $u$}
    \end{subfigure}
    \begin{subfigure}{0.4\textwidth}
        \includegraphics[width=\textwidth]{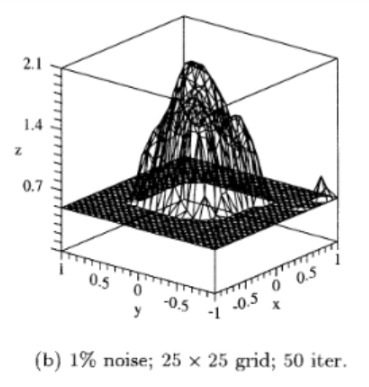}
        \caption{25 iterations, in a $25\times 25$ grid}
    \end{subfigure}

    \begin{subfigure}{0.4\textwidth}
        \includegraphics[width=\textwidth]{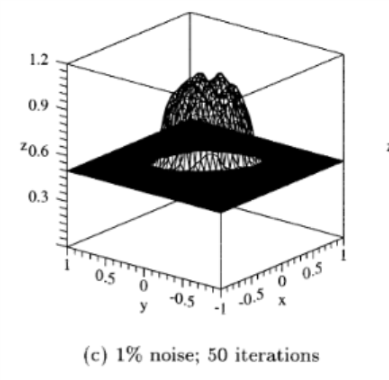}
        \caption{50 iterations, $49\times 49$ grid}
    \end{subfigure}
    \begin{subfigure}{0.4\textwidth}
        \includegraphics[width=\textwidth]{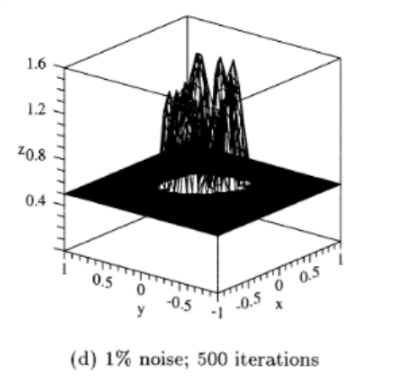}
        \caption{500 iterations, $49 \times 49$ grid}
    \end{subfigure}
    \caption{Recovery of $P$ from \cite{Knowles(PIFEP)} with $1\%$ error, for Example \ref{Example elliptic 1}.}  
    \label{knowles P recoveries}
\end{figure}
\end{Exp}{}

\begin{rem}
Note that, as mentioned previously, in this regularization method we didn't specify a lower bound for the $P$ descent, i.e., the only constraint during the minimization process is, for all $m$, $p_m|_{\partial \O} = P|_{\partial \O}$ (for the uniqueness) and not on the lower bound for $p_m$'s, where as in \cite{Knowles(PIFEP)} one has to have two constraints: (1) $p_m|_{\partial \O} = P|_{\partial \O}$ as well as $p_m \geq \nu > 0$, for some constant $\nu$, otherwise the numerical solver crashes. We can see in Figure \ref{elliptic recovery at diff iter} the initial tendency of $p_m$ towards the negative values (rather than towards the unboundedness), which is the direct manifestation of the ill-posedness in the problem. In \cite{Knowles(PIFEP)}, the remedy implemented to handle this instability is to declare a cut-off value (0.5) for the functions $p_m$, below which the values of the descent iterates are reset to the cut-off value. With this modification, the algorithm became very stable (for noise-free $u$), allowing a steady descent to the minimum, and essentially no instabilities, even after thousands iterations. However, with this new regularization method we see that stability is embedded in the process, i.e., one doesn't have to declare an external cut-off value to stabilize the process, it is self-restored (even in the case of noisy $u_\delta$).

\begin{figure}[ht]
    \centering
    \begin{subfigure}{0.4\textwidth}
        \includegraphics[width=\textwidth]{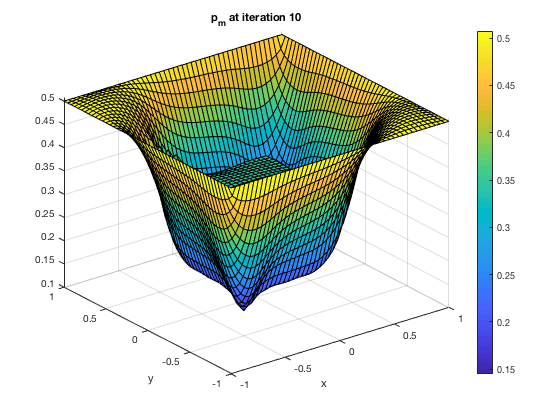}
        \caption{$p_m$ at iteration 10.}
    \end{subfigure}
    \begin{subfigure}{0.4\textwidth}
        \includegraphics[width=\textwidth]{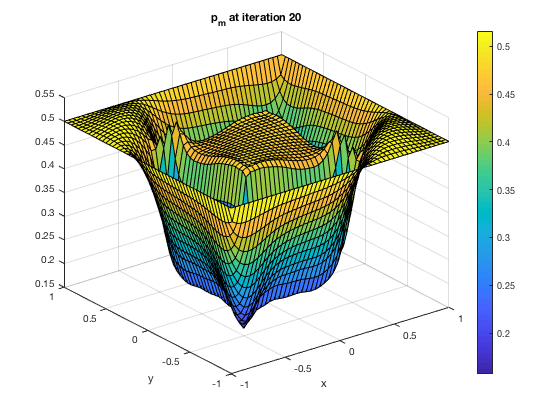}
        \caption{$p_m$ at iteration 20.}
    \end{subfigure}

    \begin{subfigure}{0.4\textwidth}
        \includegraphics[width=\textwidth]{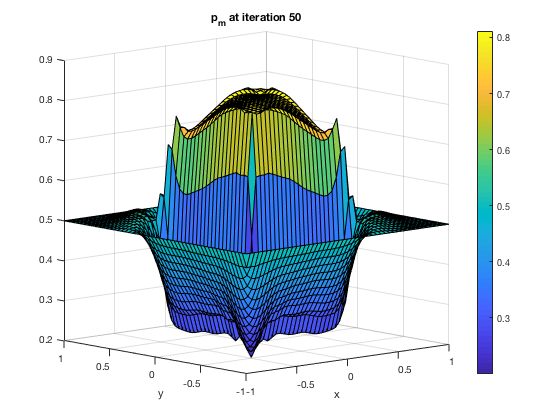}
        \caption{$p_m$ at iteration 50.}
    \end{subfigure}
    \begin{subfigure}{0.4\textwidth}
        \includegraphics[width=\textwidth]{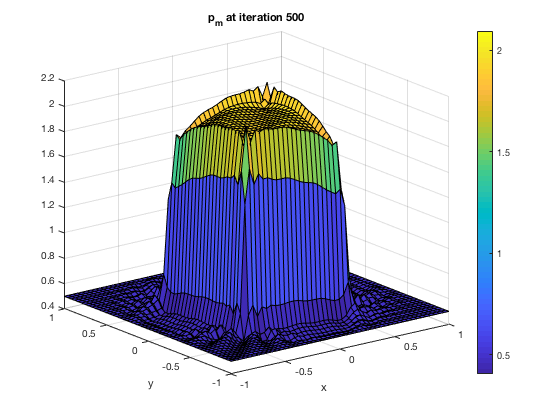}
        \caption{$p_m$ at iteration 500.}
    \end{subfigure}

    \begin{subfigure}{0.4\textwidth}
        \includegraphics[width=\textwidth]{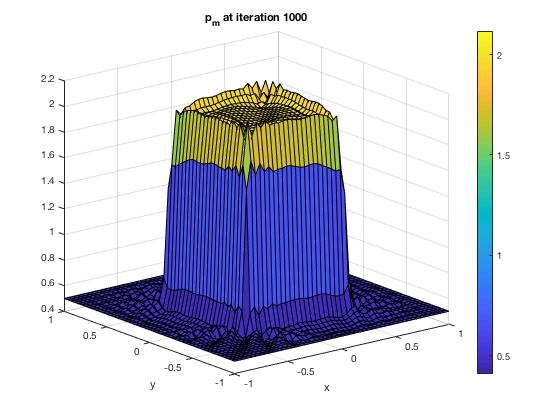}
        \caption{$p_m$ at iteration 1000.}
    \end{subfigure}
    \begin{subfigure}{0.4\textwidth}
        \includegraphics[width=\textwidth]{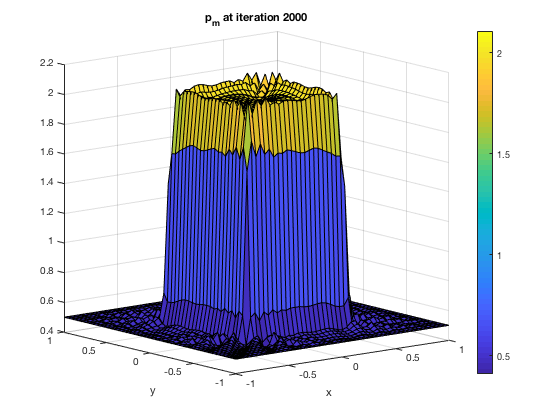}
        \caption{$p_m$ at iteration 2000.}
    \end{subfigure}
    \caption{$p_m$'s at different iterations, reflecting the stability of the process, for Example \ref{Example elliptic 1}.}  
    \label{elliptic recovery at diff iter}
\end{figure}
\end{rem}{}

The next few examples deal with the parameter $P$ having both the negative as well as the positive values, i.e., not satisfying the condition \eqref{ellipticity cond.}, but $P(x) \neq 0$ on a set of non-zero measure.

\begin{Exp}\label{Example elliptic 2}
In this example, first, we consider a no-noise situation, i.e., $u = u_\delta$ or $\delta = 0$. We consider the same domain $\O = [-1,1]\times [-1,1]$, discretized into $49\times 49$ evenly spaced grid, and the same boundary condition $\phi(x,y) = x + y + 4$ on $\partial \O$, but the parameter to be recovered here is $P(x,y) = 1 + \sin(2xy) + \cos(2xy)$. As can be seen in Figure \ref{true P2}, $P$ has both the positive and negative values (but is not zero on a set of non-zero measure) and Figure \ref{boundaryP2} shows (the given) $P|_{\partial \O}$ such that $p_m|_{\partial \O} = P|_{\partial \O}$, for all $m$. We don't contaminate the data $u$ with any noise and perform the minimization. The recovered parameter is shown in Figure \ref{boundaryP2}, with a relative error of 0.0033\%.

\begin{figure}[ht]
    \centering
    \begin{subfigure}{0.4\textwidth}
        \includegraphics[width=\textwidth]{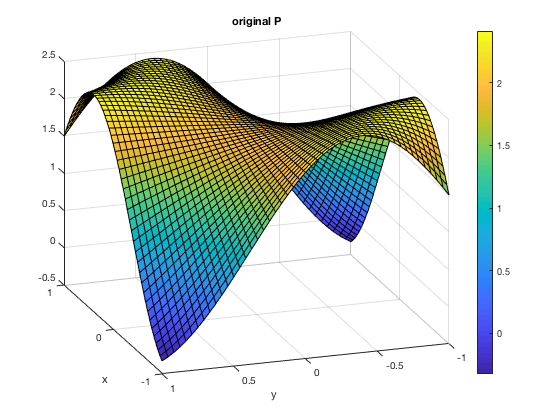}
        \caption{True $P$}
        \label{true P2}
    \end{subfigure}
    \begin{subfigure}{0.4\textwidth}
        \includegraphics[width=\textwidth]{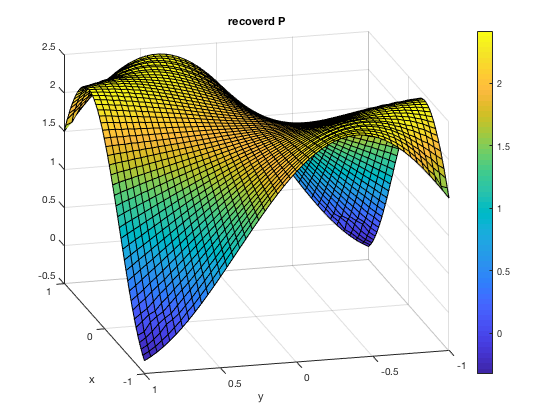}
        \caption{Recovered $\tilde{P}$}
        \label{recovered P2}
    \end{subfigure}
    \caption{True $P$ and recovered $\tilde{P}$, for Example \ref{Example elliptic 2}.}  
    \label{true and recovered P2}
\end{figure}

We repeat this example but this time with a noisy $u_\delta$, such that $\frac{\Lonenorm{u_\delta - u}}{\Lonenorm{u}} \approx 0.74\%$. The reason for such small error is due to the presence of two spikes in computed $u$, see Figure \ref{noisy U2}, which is a result of the presence of positive and negative values in $P$ and hence, a large error would have lead to losing of the spikes. Also note that in this case one can not use a polynomial surface-fit, as it will again lose the spikes. In this case we simply fit a surface ($\tilde{u_\delta}$), see Figure \ref{smoothed U2}, generated by piecewise cubic interpolation, the disadvantage being, this leads to huge errors when estimating the $\nabla u$ using $\nabla \tilde{u_\delta}$ and also, one also has to be very careful when computing $T_{u_\delta}(p_m) = -\nabla \cdot (p_m \nabla \tilde{u_\delta})$. Anyway, after carefully handling all of the above concerns, the recovered $P$ is shown in Figure \ref{recovered P2_1}, where the relative error in the recovery is around 25.57\%.   

\begin{figure}[ht]
    \centering
    \begin{subfigure}{0.4\textwidth}
        \includegraphics[width=\textwidth]{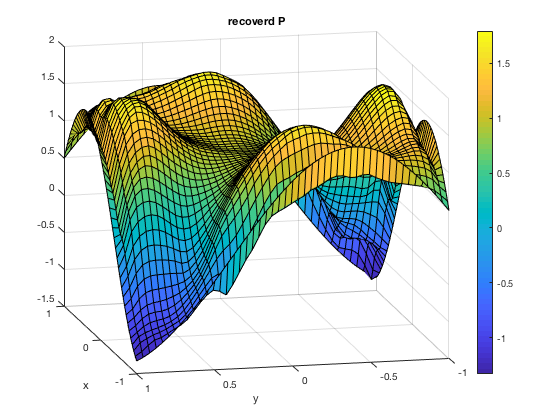}
        \caption{Recovered $\tilde{P}$, for Example \ref{Example elliptic 2}.}
        \label{recovered P2_1}
    \end{subfigure}
    \begin{subfigure}{0.4\textwidth}
        \includegraphics[width=\textwidth]{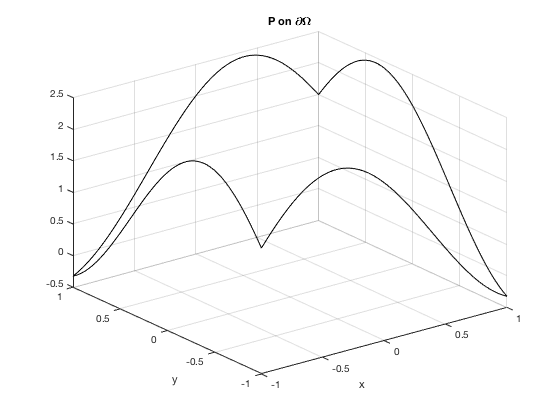}
        \caption{$P$ on $\partial \O$ for Example \ref{Example elliptic 2}.}
        \label{boundaryP2}
    \end{subfigure}
    \caption{$P$ on the boundary $\partial \O$.}  
    \label{recoveredP2_1 and boundary P2}
\end{figure}

\begin{figure}[ht]
    \centering
    \begin{subfigure}{0.4\textwidth}
        \includegraphics[width=\textwidth]{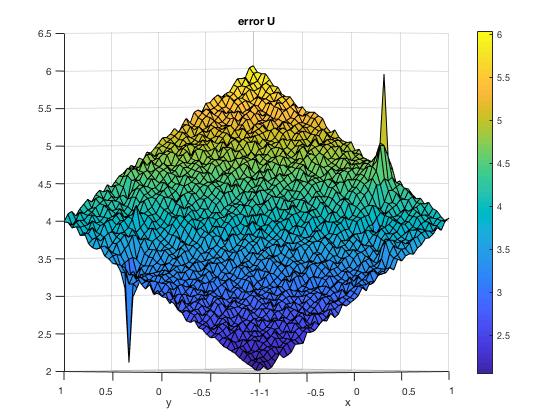}
        \caption{noisy $u_\delta$.}
        \label{noisy U2}
    \end{subfigure}
    \begin{subfigure}{0.4\textwidth}
        \includegraphics[width=\textwidth]{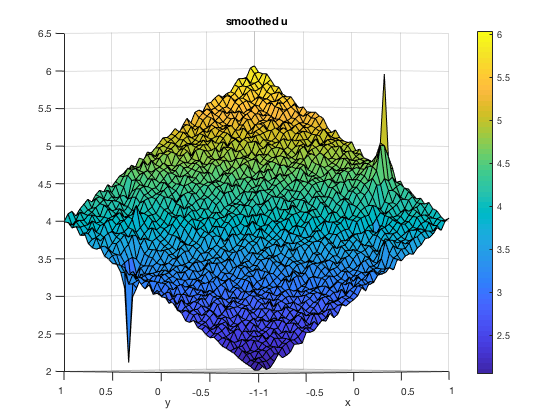}
        \caption{smoothed $\tilde{u_\delta}$}
        \label{smoothed U2}
    \end{subfigure}
    \caption{Noisy and smoothed data for Example \ref{Example elliptic 2}}  
    \label{noisy and smoothed U2}
\end{figure}
\end{Exp}{}

\begin{Exp}\label{Example elliptic 3}
Here we consider another example of positive and negative $P$, which is defined on $\O$ as
\beq
    P(x,y) = 
    \begin{cases}
    -2, \;\; &\mbox{ if } |x| < 0 \;\; \mbox{ and } |y| < 0\\
    0.5, \;\; &\mbox{ if } |x| \geq 0 \;\; \mbox{ and } |y| < 0\\
    0.5, \;\; &\mbox{ if } |x| < 0 \;\; \mbox{ and } |y| \geq 0\\
    2, \;\; &\mbox{ if } |x| \geq 0 \;\; \mbox{ and } |y| \geq 0. 
    \end{cases}{}
\eeq
The choice of the above $P$ increased the complexity of the computed solution in manifold ways, see Figure \ref{noerrorU3}. Hence, we do not add any external noise in this example, i.e., $u_\delta = u$. The recovered $P$, with a relative error of 40\%, is shown in Figure \ref{true and recovered P3}, and the constraint $p_m|_{\partial \O} = P|_{\partial \O}$ for the minimization is shown in Figure \ref{boundaryP3}

\begin{figure}[ht]
    \centering
    \begin{subfigure}{0.4\textwidth}
        \includegraphics[width=\textwidth]{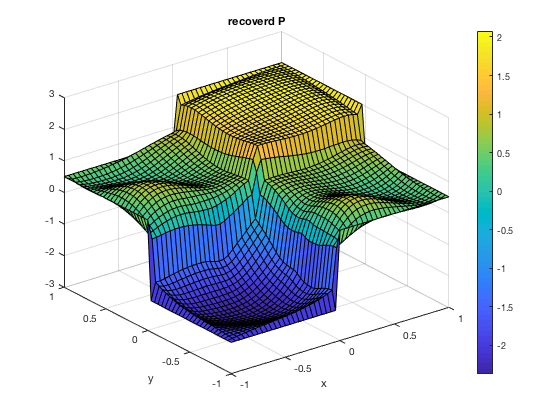}
        \caption{Recovered $\tilde{P}$, for Example \ref{Example elliptic 3}.}
        \label{recovered P3}
    \end{subfigure}
    \begin{subfigure}{0.4\textwidth}
        \includegraphics[width=\textwidth]{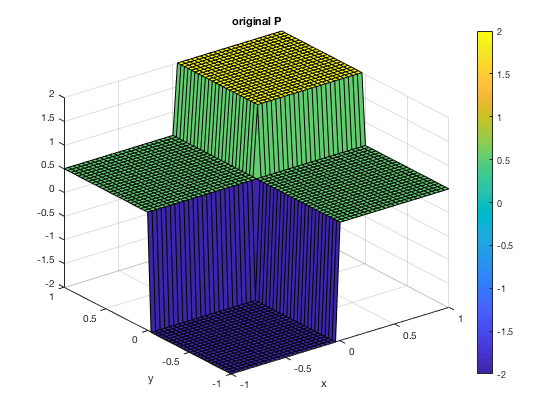}
        \caption{True $P$ for Example \ref{Example elliptic 3}.}
        \label{trueP3}
    \end{subfigure}
    \caption{True and recovered $P$ for Example \ref{Example elliptic 3}.}  
    \label{true and recovered P3}
\end{figure}

\begin{figure}[ht]
    \centering
    \begin{subfigure}{0.4\textwidth}
        \includegraphics[width=\textwidth]{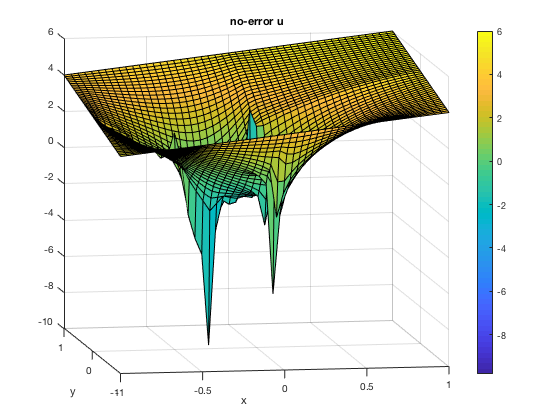}
        \caption{No-error $u$ for Example \ref{Example elliptic 3}.}
        \label{noerrorU3}
    \end{subfigure}
    \begin{subfigure}{0.4\textwidth}
        \includegraphics[width=\textwidth]{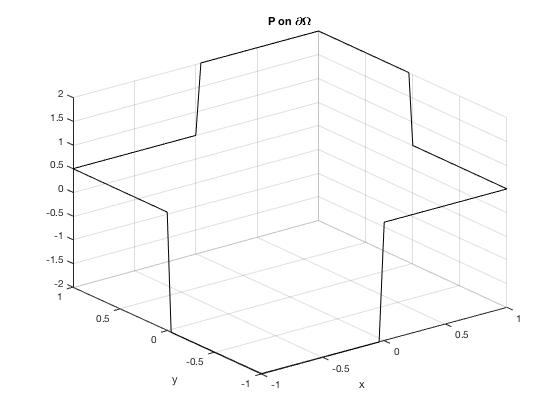}
        \caption{$P$ on $\partial \O$ for Example \ref{Example elliptic 3}.}
        \label{boundaryP3}
    \end{subfigure}
    \caption{ \hspace{0.1cm} $u$ and $P|_{\partial \O}$ for Example \ref{Example elliptic 3}.}  
    \label{noerrorU3 and boundaryP3}
\end{figure}
\end{Exp}{}

\begin{Exp}\label{Example elliptic 4}
In the final example we demonstrate the importance of the boundary information. As seen in the Figures \ref{boundaryP2} and \ref{boundaryP3}, the boundary data $P|_{\partial \O}$ does provide us the information that the parameter $P$ has both the positive and negative values. In this example we choose a $P$ which has both the negative and positive values, but the boundary data $P|_{\partial \O}$ doesn't reflect it, i.e., $P$ is positive and negative in a strictly interior region $\O' \subsetneq \O \backslash \partial \O$. The parameter $P$ is defined as 
\beq
    P(x,y) = 
    \begin{cases}
    -2, \;\; &\mbox{ if } -0.25 < |x| < 0.75 \;\; \mbox{ and } -0.25 < |y| < 0.75\\
    2, \;\; &\mbox{ if } 0.25 < |x| < 0.75 \;\; \mbox{ and } 0.25 < |y| < 0.75\\
    1, \;\; &\mbox{ otherwise}.
    \end{cases}{}
\eeq
So one can observe that, from Figures \ref{trueP4} or \ref{boundaryP4}, $P|_{\partial \O}$ does not contain any information regarding the negativity of $P$. Again, since the computed $u$, see Figure \ref{noerrorU4}, has a complex structure we do not impose additional noise to it. Figure \ref{recovered P4} shows the recovered $\tilde{P}$, with a relative error of 65.22\% (after 704 iterations), and Figure \ref{boundaryP4} shows boundary data $P|_{\partial \O}$. Note that, as mentioned above, the recovery is unstable for $p_m = 0$ on a set of non-zero measure and hence, the descent process is extremely slow when $p_m$ approaches 0 from the positive side, in an attempt to cross over to the negative side; where as in the previous examples (Example \ref{Example elliptic 2} and \ref{Example elliptic 3}), since $p_m|_{\partial \O} = P|_{\partial \O}$ has both the positive and negative values, $p_m|_{p_m > 0 } \rightarrow P|_{P > 0}$ and $p_m|_{p_m < 0} \rightarrow P|_{P < 0}$ (weakly) in $\L$.

\begin{figure}[ht]
    \centering
    \begin{subfigure}{0.4\textwidth}
        \includegraphics[width=\textwidth]{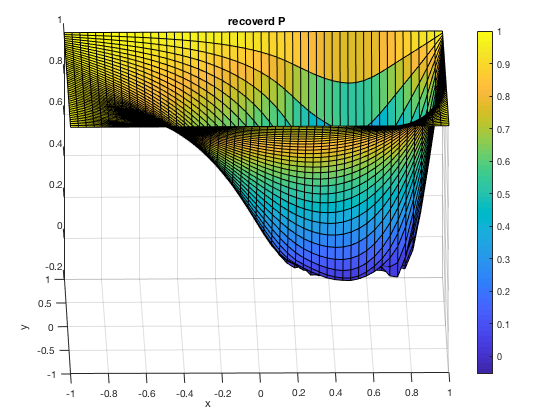}
        \caption{Recovered $\tilde{P}$, for Example \ref{Example elliptic 4}.}
        \label{recovered P4}
    \end{subfigure}
    \begin{subfigure}{0.4\textwidth}
        \includegraphics[width=\textwidth]{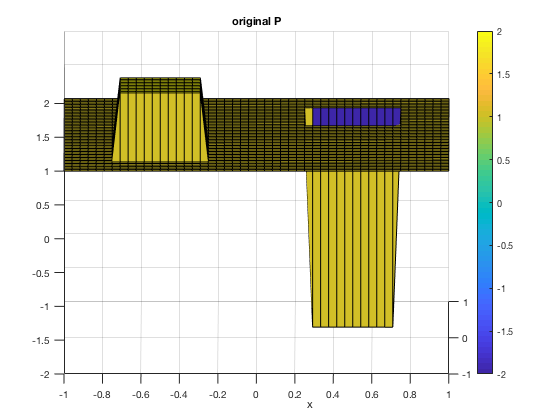}
        \caption{True $P$ for Example \ref{Example elliptic 4}.}
        \label{trueP4}
    \end{subfigure}
    \caption{\hspace{0.1cm} True and recovered $P$ for Example \ref{Example elliptic 4}.}  \label{true and recovered P4}
\end{figure}

\begin{figure}[ht]
    \centering
    \begin{subfigure}{0.4\textwidth}
        \includegraphics[width=\textwidth]{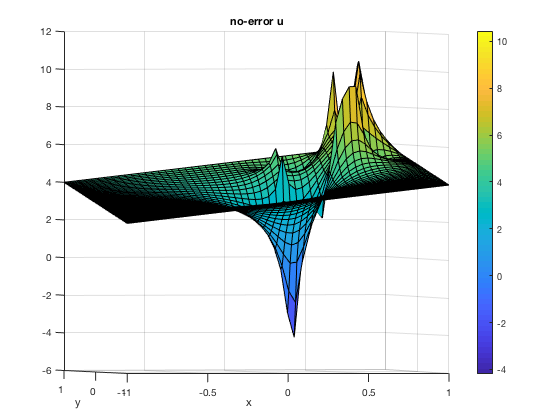}
        \caption{No-error $u$ for Example \ref{Example elliptic 4}.}
        \label{noerrorU4}
    \end{subfigure}
    \begin{subfigure}{0.4\textwidth}
        \includegraphics[width=\textwidth]{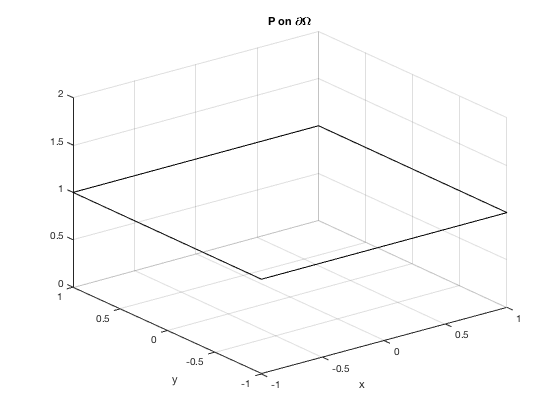}
        \caption{$P$ on $\partial \O$ for Example \ref{Example elliptic 4}.}
        \label{boundaryP4}
    \end{subfigure}
    \caption{\hspace{0.1cm} $u$ and $P|_{\partial \O}$ for Example \ref{Example elliptic 4}.}  
    \label{noerrorU4 and boundaryP4}
\end{figure}
\end{Exp}{}

\begin{rem}
As we can see that, using the developed regularization method, the recovery of the parameter $P$ is very efficient even in the presence of extreme errors or when $P$ does not satisfy \eqref{ellipticity cond.}. Though we were not able to efficiently estimate $\nabla u_\delta$, still the method was very stable and quite effective. Hence, we expect to have even better results if we can extend the numerical differentiation procedure, developed for the single variable in Chapter $??$, to the multi-variable scenario. The other area of interest is to handle the case when the boundary data $P|_{\partial \O}$ does not contain any information about the parameter $P$ having both the positive and negative values. One may make use of the solution $u$ in that case, since the solution has peaks if $P$ has both positive and negative values. 
\end{rem}{}

\bibliography{thesisref} 
\bibliographystyle{ieeetr}

\end{document}